\documentclass[12pt]{amsart}
\usepackage{amsfonts}
\usepackage{amssymb}
\numberwithin{equation}{section}
\theoremstyle{plain}
 \newtheorem{thm}{Theorem}[section]
 \newtheorem{lem}{Lemma}[section]
 \newtheorem{cor}{Corollary}[section] 
 
\theoremstyle{definition}
 \newtheorem{defn}{Definition}[section]
 
 \newtheorem{rem}{Remark}[section]

\newcommand{\ep}{\varepsilon}

\newcommand{\mcal}{\mathcal}

\newcommand{\q}{\quad}

\newcommand{\R}{\mathbb{R}}
\newcommand{\N}{\mathbb{N}}

\newcommand{\M}{\mathbf{M}}

\begin{document}
\setlength{\baselineskip}{15pt}
\setlength{\parindent}{1.8pc}

\title{The conjectures of Embrechts and Goldie}
\author{Toshiro Watanabe }
\maketitle
\footnote[0]{T. Watanabe: Center for Mathematical Sciences The Univ. of Aizu. Aizu-Wakamatsu 965-8580, Japan. e-mail: t-watanb@u-aizu.ac.jp\\ }

{\small
{\bf Abstract}.  It is shown that the class of convolution equivalent distributions and the class of locally  subexponential distributions are not closed under convolution roots. Moreover, two sufficient conditions for  the closure under convolution roots of the class of convolution equivalent distributions are given. }\\

\medskip  

             

\section{Introduction and main results}
In what follows, we denote by  $\R$ the real line and by  $\R_+$ the half line $[0,\infty)$.  Let $\N$ be the totality of positive integers. 
Let $\eta$ and $\rho$ be probability measures on $\R$. We denote the convolution of $\eta$ and $\rho$ by $\eta*\rho$ and denote the $n$-th convolution power of $\rho$ by  $\rho^{n*}$.  
 Let $f(x)$ and $g(x)$ be integrable functions on $\R$. We denote the convolution of $f(x)$ and $g(x)$ by $f*g(x)$ and  denote the $n$-th convolution power of  $f(x)$ by $f^{n*}(x)$.  For positive  functions $f_1(x)$ and $g_1(x)$ on $[a,\infty)$ for some $a \in \R$, we define the relation $f_1(x) \sim g_1(x)$ by $\lim_{x \to \infty}f_1(x)/g_1(x) =1$ and the relation $f_1(x) \asymp g_1(x)$ by $0 <\liminf_{x \to \infty}f_1(x)/g_1(x) \le\limsup_{x \to \infty}f_1(x)/g_1(x) < \infty$.   The  tail of a probability  measure $\eta$ on $\R$ is denoted by  $\bar \eta(x)$, that is, $\bar \eta(x): = \eta((x,\infty))$ for $x \in \R$. Let $\gamma \in \R$. The $\gamma$-exponential moment of $\eta$ is denoted by $\widehat \eta(\gamma)$, namely, 
$$\widehat \eta(\gamma):= \int_{-\infty}^{\infty}e^{\gamma x}\eta(dx).$$ 
If $\widehat \eta(\gamma) < \infty$, we define for $z \in \R$
$$\widehat \eta(\gamma +iz):= \int_{-\infty}^{\infty}e^{(\gamma +iz)x}\eta(dx).$$ 
We define the class  ${\mcal P}_+$ as the totality of probability distributions on $\R_+$.  We use the words "increase" and "decrease" in the wide sense allowing flatness.

\begin{defn}\quad Let $\gamma \ge 0$.

(i) A distribution $\rho$ on  $\R$ belongs to the class $\mcal{ L}(\gamma)$ if $\bar \rho(x) >0$ for every $x \in  \R$ and if
\begin{equation}\label{2.1}
 \bar \rho(x+a)\sim  e^{-\gamma a}\bar \rho(x) \q \mbox{ for every  }\q  a \in \R. \nonumber
\end{equation}

(ii)  A distribution $\rho$ on $\mathbb{R}$ belongs to the class ${\mcal S}(\gamma)$ if $\rho \in {\mcal L}(\gamma)$ with $\widehat \rho(\gamma)< \infty$ and if
\begin{equation}\label{2.2}
 \overline{\rho^{2*}}(x) \sim  2\widehat \rho(\gamma)\bar \rho(x). \nonumber
\end{equation}

(iii) Let $\gamma_1 \in \R$. A distribution $\rho$ on $\mathbb{R}$ belongs to the class ${\mcal M}(\gamma_1)$ if $\widehat \rho(\gamma_1)< \infty$. 

\end{defn}
\medskip

\begin{defn}

(i)  A nonnegative measurable function  $g(x)$  on $\R$  belongs to the class $ {\bf L}$   if $g(x)>0$ for all sufficiently large $x>0$ and if $g(x+a) \sim g(x)$ for any $a \in \R$.

(ii) A nonnegative measurable function  $g(x)$  on $\R$  belongs to the class ${\bf M}$   if it is positive for  all sufficiently large $x>0$ and there exist $a \in [0,1]$, $f_1(x) \in  {\bf L}$, and a positive and  decreasing function $f_2(x)$ on $\R$ such that  $g(x) \sim af_1(x)+ (1-a)f_2(x)$.
\end{defn}
\begin{defn}

(i)   Let $\Delta:= (0,c]$ with $c>0.$ A distribution $\rho$   on  $\R$   belongs to the class $ {\mcal L}_{\Delta}$   if $\rho((x,x+c])  \in   {\bf L}$.
 
(ii)  Let $\Delta:= (0,c]$ with $c>0.$  A distribution $\rho$   on  $\R$   belongs to the class $ {\mcal S}_{\Delta}$   if $\rho  \in  {\mcal L}_{\Delta}$ and 
$\rho^{2*}((x,x+c])  \sim 2 \rho((x,x+c]). $ 

(iii) A distribution  $\rho$  on  $\R$    belongs to the class ${\mcal L}_{loc}$   if $\rho  \in {\mcal L}_{\Delta}$ for each  $\Delta:= (0,c]$ with $c>0.$ 

(iv)  A distribution  $\rho$   on  $\R$  belongs to the class ${\mcal S}_{loc}$  if $\rho  \in {\mcal S}_{\Delta}$ for each $\Delta:= (0,c]$ with $c>0.$

\end{defn}
Distributions in the class ${\mcal S}(\gamma)$ are called {\it convolution equivalent} and those in the class ${\mcal S}_{loc}$ are called {\it locally subexponential}.  The study of the class ${\mcal S}(\gamma)$ goes back to Chover et al.\ \cite{cnw,cnw1}. The class $ {\mcal S}_{\Delta}$ is introduced by Asmussen et al.\ \cite{afk03} and the class ${\mcal S}_{loc}$  is by Watanabe and Yamamuro \cite{wy10}. Applications of those classes include renewal theory, random walks, queues, branching processes, L\'evy processes, and infinite divisibility.

\begin{defn}
 We say that a class $\mathcal{C}$ of probability distributions  on $\R$ is closed under convolution roots if  $\mu^{n*} \in \mathcal{C}$ for some $ n \in \mathbb{N}$ implies that $\mu \in \mathcal{C}$. 

\end{defn}
 Embrechts et al.\ \cite{egv79} in the one-sided case and   Watanabe \cite{wa08} in the two-sided case proved that the class $\mathcal{S}(0)$ of {\it subexponential distributions} is closed under convolution roots. Embrechts and Goldie stated in \cite{eg82} that a crucial point for proving limit theorems using  $\mathcal{S}(\gamma)$ is the convolution roots closure of $\mathcal{S}(\gamma)$. Further,  they gave the following conjectures in \cite{eg80,eg82}, respectively.

\vspace{1ex}

{\it Conjecture I. } The class ${\mcal L}(\gamma)$  with $\gamma \ge 0$    is closed under convolution roots. 

\vspace{1ex}

{\it Conjecture II. } The class $\mathcal{S}(\gamma)$ with $\gamma > 0$  is  closed under convolution roots.

\vspace{1ex}

 Embrechts and Goldie \cite{eg82} in the one-sided case and Pakes \cite{pa07} in the two-sided case  obtained the following.

\vspace{1ex}

{\bf Theorem A} {\it  Let $\gamma >0$ and let $\mu$ be a distribution   on $\R$. If $\mu \in \mathcal{L}(\gamma)$ and $\mu^{n*} \in \mathcal{S}(\gamma)$ for some $ n \in \mathbb{N}$, then $\mu \in \mathcal{S}(\gamma)$}.

\vspace{1ex}
Moreover,  Watanabe showed in Theorem 1.1 of \cite{wa08} the following.

\vspace{1ex}

{\bf Theorem B} {\it  Let $\gamma >0$ and let $\mu$ be an   infinitely divisible distribution on $\R$. If $\mu^{n*} \in \mathcal{S}(\gamma)$ for some $ n \in \mathbb{N}$, then $\mu \in \mathcal{S}(\gamma)$.}

\vspace{1ex}

We see from Theorem A that if  Conjecture I is true for every $\gamma > 0$, then so is Conjecture II.  However,  Shimura and Watanabe \cite{shw} disproved  Conjecture I for every $\gamma \ge 0$. On the other hand, Conjecture II was unsolved for over thirty years. In this paper, we disprove  Conjecture II for every $\gamma > 0$.  Classical Wiener's approximation theorem plays a key role for the solution. The closure problem under  convolution roots for  the other distribution classes is discussed by Shimura and Watanabe \cite{shw05} and  Watanabe and Yamamuro  \cite{wy09,wy10a}.  Our main results are as follows.
\begin{thm}
 The class $\mathcal{S}(\gamma)$ with $\gamma > 0$ is not closed under convolution roots.
\end{thm} 
  By using Lemma 2.4 below, we have the following corollary.
\begin{cor}  Let $\Delta:= (0,c]$ with $c>0$ and let $\gamma > 0$. We have the following.

(i) The class $\mathcal{S}_{loc}$ is not closed under convolution roots.

(ii)  The class $\mathcal{S}_{\Delta}$ is not closed under convolution roots.

(iii)  The class $\mathcal{L}(\gamma)\cap \mathcal{M}(\gamma)$ is not closed under convolution roots.

(iv) The class $\mathcal{L}_{loc}$ is not closed under convolution roots.

(v)  The class $\mathcal{L}_{\Delta}$ is not closed under convolution roots.
\end{cor}

Next,  we establish an extension of Theorem A. Note that the class ${\bf L}$ is a proper subclass of the class ${\bf M}$ and that  $\mu \in \mathcal{L}(\gamma)$ with $\gamma > 0$ if and only if $e^{\gamma x}\bar\mu(x) \in {\bf L}$.

\begin{thm} Let $\gamma >0$ and let $\mu$ be a distribution   on $\R$. Assume that  $e^{\gamma x}\bar\mu(x) \in \M$. Then, $\mu^{n*} \in \mathcal{S}(\gamma)$  for some $n \in \N$ implies that  $\mu \in \mathcal{S}(\gamma)$. 

\end{thm}
\begin{cor}
 Let $\mu$ be a distribution   on $\R$. Assume that $\widehat \mu(-\gamma)< \infty$ for some $\gamma >0$ and that $\mu((x,x+c]) \in \M$ for every $c >0$. Then, $\mu^{n*} \in \mathcal{S}_{loc}$  for some $n \in \N$ implies that  $\mu \in \mathcal{S}_{loc}$. 
\end{cor}

Finally, we present an extension of Theorem B. Note that if $\widehat\mu(\gamma ) < \infty$ for an infinitely divisible distribution $\mu$ on $\R$, then $\widehat\mu(\gamma + iz) \neq 0$ for every $z \in \R$. See Theorem 25.17  of Sato \cite{sato}.
\begin{thm} Let $\gamma >0$ and let $\mu$ be a distribution   on $\R$.
Assume that $\widehat\mu(\gamma ) < \infty$ and $\widehat\mu(\gamma + iz) \neq 0$ for every $z \in \R$. Then, $\mu^{n*} \in \mathcal{S}(\gamma)$ for some $ n \in \mathbb{N}$ implies that $\mu \in \mathcal{S}(\gamma)$. 
\end{thm}

\begin{rem} We see from Theorem 1.3 that the distribution $\mu$ on $\R$ of each counter-example to Conjecture II must satisfy that $\widehat\mu(\gamma ) < \infty$ and $\widehat\mu(\gamma + iz_0) = 0$ for some $z_0 \in \R$.
\end{rem}

\begin{cor} Let $\mu$ be a distribution   on $\R$.
Assume that $\widehat\mu(-\gamma ) < \infty$ for some $\gamma >0$ and $\widehat\mu(iz) \neq 0$ for every $z \in \R$. Then,  $\mu^{n*} \in \mathcal{S}_{loc}$ for some $ n \in \mathbb{N}$ implies that $\mu \in \mathcal{S}_{loc}$.
\end{cor}

In Sect.\ 2, we give preliminaries for the proofs of the main results.   In Sect.\ 3, 4, and 5, we prove Theorems 1.1, 1.2, and 1.3 and their corollaries, respectively.  

\section{Preliminaries}
We define the class ${\mcal L}_{ac}$ of distributions with {\it long-tailed densities} and  the class ${\mcal S}_{ac}$ of distributions with {\it subexponential densities}.
\begin{defn}

(i)  A distribution $\rho$  on $\R$   belongs to the class ${\mcal L}_{ac}$ if there is $g(x)\in {\bf L}$ such that $\rho(dx) =g(x) dx$.

(ii) A distribution $\rho$   on $\R$   belongs to the class ${\mcal S}_{ac}$  if  there is $g(x)\in {\bf L}$ such that  $\rho(dx) =g(x) dx$ and
$g^{2*}(x)  \sim 2 g(x)$.

\end{defn}
The class ${\mcal S}_{ac}$ (resp. ${\mcal L}_{ac}$) is a proper subclass of the class  ${\mcal S}_{loc}$ (resp. ${\mcal L}_{loc}$).

\begin{lem} 

(i) Let $f(x)dx$ be a distribution on $\R_+$. If $f(x) \sim c x^{-\alpha}$ with $c >0$ and $\alpha >1$, then  $f(x)dx \in {\mcal S}_{ac}$.

(ii) Let $\gamma >0$ and $\mu$ be a distribution on $\R_+$. If $\bar \mu(x) \sim ce^{-\gamma x}x^{-\alpha}$ with $c >0$ and $\alpha >1$, then  $\mu \in {\mcal S}(\gamma)$.

\end{lem}

{\it Proof } The proof of assertion (i) is clear from Theorem 4.14 of Foss et al.\ \cite{fkz13}. The proof of assertion (ii) is due to  assertion (i) and Theorem 2.1  of Kl\"uppelberg \cite{kl89}. \qed

\begin{lem} Let $\gamma \ge 0$. We have the following.

(i) Let $\mu \in \mathcal{L}(\gamma)$ with ${\widehat \mu}(\gamma) < \infty$.  Then, $\mu \in \mathcal{S}(\gamma)$ if and only if 
\begin{equation}
\lim_{A \to \infty}\limsup_{x \to \infty}\frac{\int_{A+}^{(x-A)+}{\overline{\mu}(x-u)}\mu(du)}{\overline{\mu}(x)}=0.\nonumber
\end{equation}

(ii) Let $\mu_1$ and $\mu_2$ be distributions on $\R$. If $\mu_1 \in \mathcal{S}(\gamma)$ and $\overline{\mu_2}(x) \sim c \overline{\mu_1}(x)$ with $c >0$, then 
$\mu_2 \in \mathcal{S}(\gamma)$. 

\end{lem}
{\it Proof } \q  First, we prove assertion (i). Let $\mu \in \mathcal{L}(\gamma)$ with ${\widehat \mu}(\gamma) < \infty$ and let $A>0$. We have, for $x > 2A$,
\begin{equation}
\overline{\mu^{2*}}(x) = \sum_{j=1}^3 H_j(x), \nonumber
\end{equation}
where 
\begin{equation}
H_1(x) := 2\int_{-\infty}^{A+}\overline{ \mu}(x-u)\mu(du), \q H_2(x) :=\overline{ \mu}(x-A)\overline{ \mu}(A),\nonumber
\end{equation}
and
\begin{equation}
H_3(x) := \int_{A+}^{(x-A)+}\overline{ \mu}(x-u)\mu(du).\nonumber
\end{equation}
Since we see that
\begin{equation}
\sup_{u \in (-\infty,A]}\frac{\overline{ \mu}(x-u)}{\overline{ \mu}(x)}\le \frac{\overline{ \mu}(x-A)}{\overline{ \mu}(x)} \mbox{ and } \lim_{x \to \infty}\frac{\overline{ \mu}(x-A)}{\overline{ \mu}(x)}= e^{\gamma A},\nonumber
\end{equation}
we obtain from the dominated convergence theorem that
\begin{equation}
\lim_{x \to \infty} \frac{H_1(x) }{\overline{ \mu}(x)}=2\int_{- \infty}^{A+}\lim_{x \to \infty}\frac{\overline{ \mu}(x-u)}{\overline{ \mu}(x)}\mu(du)= 2\int_{- \infty}^{A+}e^{\gamma u}\mu(du).
\end{equation}
We have 
\begin{equation}
 \lim_{A \to \infty}\lim_{x \to \infty} \frac{H_2(x) }{\overline{ \mu}(x)}= \lim_{A \to \infty}e^{\gamma A}\overline{ \mu}(A) \le  \lim_{A \to \infty}\int_{A+}^{\infty}e^{\gamma x}\mu(dx)=0.
\end{equation}
Thus, we see from (2.1) and (2.2) that
\begin{equation}
\begin{split}
0 = \lim_{A \to \infty}\lim_{x \to \infty} \frac{H_1(x) }{\overline{ \mu}(x)}-2{\widehat \mu}(\gamma) & \le \liminf_{x \to \infty}\frac{\overline{\mu^{2*}}(x)}{\overline{ \mu}(x)}-2{\widehat \mu}(\gamma)\\ \nonumber
 & \le \limsup_{x \to \infty}\frac{\overline{\mu^{2*}}(x)}{\overline{ \mu}(x)}-2{\widehat \mu}(\gamma)\\ \nonumber
 & =\lim_{A \to \infty} \limsup_{x \to \infty}\frac{H_3(x) }{\overline{ \mu}(x)}.\nonumber
\end{split}
\end{equation}
Thus, we find that $\mu \in \mathcal{S}(\gamma)$, that is,  
$$\lim_{x \to \infty}\frac{\overline{\mu^{2*}}(x)}{\overline{ \mu}(x)}=2{\widehat \mu}(\gamma)$$ 
if and only if
 $$ \lim_{A \to \infty} \limsup_{x \to \infty}\frac{H_3(x) }{\overline{ \mu}(x)}=0.$$ 
The proof of assertion (ii) is due to Lemma 2.4 of Pakes \cite{pa04}. \qed

\medskip

Let $\gamma \in \R$. For $ \mu \in  {\mcal M}(\gamma)$, we define the {\it exponential tilt} $\mu_{\langle\gamma\rangle}$ of $\mu$ as

\begin{equation}
\mu_{\langle\gamma\rangle}(dx):= \frac{1}{\widehat\mu(\gamma)}e^{\gamma x}\mu(dx). \nonumber
\end{equation}
Exponential tilts preserve convolutions, that is, $(\mu*\rho)_{\langle\gamma\rangle}=\mu_{\langle\gamma\rangle}*\rho_{\langle\gamma\rangle}$ for distributions $\mu, \rho \in {\mcal M}(\gamma)$. Let  ${\mcal C}$ be a distribution class. For a  class ${\mcal C}\subset {\mcal M}(\gamma)$, we define the class ${\frak E}_{\gamma}( {\mcal C})$ by
\begin{equation}
{\frak E}_{\gamma}( {\mcal C}):=\{\mu_{\langle\gamma\rangle} : \mu \in {\mcal C}\}. \nonumber
\end{equation}
It is obvious that ${\frak E}_{\gamma}(  {\mcal M}(\gamma))=  {\mcal M}(-\gamma) $ and that $(\mu_{\langle\gamma\rangle})_{\langle-\gamma\rangle}= \mu$ for $\mu \in {\mcal M}(\gamma)$.
 The class ${\frak E}_{\gamma}( {\mcal S}(\gamma))$ is determined  by Watanabe and Yamamuro as follows.

\begin{lem} (Theorem 2.1 of \cite{wy10})
Let $\gamma >0$.  We have the following.

(i) We have ${\frak E}_{\gamma}( {\mcal L}(\gamma)\cap {\mcal M}(\gamma))= {\mcal L}_{loc}\cap {\mcal M}(-\gamma) $ and hence ${\frak E}_{\gamma}( {\mcal L}(\gamma)\cap {\mcal M}(\gamma)\cap{\mcal P}_+)= {\mcal L}_{loc}\cap{\mcal P}_+ $. Moreover, if $ \rho \in {\mcal L}(\gamma)\cap {\mcal M}(\gamma)$, then we have 
\begin{equation}
\rho_{\langle\gamma\rangle}((x,x+c]) \sim \frac{c\gamma}{\widehat\rho(\gamma)}e^{\gamma x}\bar\rho(x) \mbox{ for all } c>0. \nonumber
\end{equation}

(ii)  We have ${\frak E}_{\gamma}( {\mcal S}(\gamma))= {\mcal S}_{loc}\cap {\mcal M}(-\gamma) $ and thereby  ${\frak E}_{\gamma}( {\mcal S}(\gamma)\cap{\mcal P}_+)= {\mcal S}_{loc}\cap{\mcal P}_+ $.
\end{lem} 
 
\medskip

A straightforward consequence of the above lemma is  the following.

\begin{lem} Let $\gamma >0$. We have the following.

(i) The class $\mathcal{S}(\gamma)$ is closed under convolution roots if and only if so is the class  ${\mcal S}_{loc}\cap {\mcal M}(-\gamma) $.

(ii) The class $\mathcal{L}(\gamma)\cap {\mcal M}(\gamma) $ is closed under convolution roots if and only if so is the class  ${\mcal L}_{loc}\cap {\mcal M}(-\gamma) $.

\end{lem}

In Sect.\ 5, we shall use the following lemma, namely, Wiener's approximation theorem in \cite{w32} for the proof of Theorem 1.3.

\begin{lem} (Theorem 4.8.4 of  \cite{bi87} or Theorem 8.1 of  \cite{ko04}) For $f(x) \in L^1(\R),$ the following are equivalent:

(1) $\int_{-\infty}^{\infty}\exp(izx)f(x)dx \ne 0 \mbox{ for every } z \in \R.$

(2) If, for a bounded measurable function $g(x)$ on $\R$,
$$ \int_{-\infty}^{\infty}g(x-t)f(t)dt =0  \mbox{ for every } x \in \R,$$ 
then $g(x)=0$ for almost every $x \in \R$.
\end{lem}

\section{Proofs of Theorem 1.1  and its corollary} 
We prove Theorem 1.1 only for $\gamma =1$. The general case for  $\gamma >0$ is similar and omitted. The symbol $[x]$ stands for the largest integer not exceeding a real number $x$ and the symbol $1_{B}(x)$ does for the indicator function of a subset $B$ of $\R$.  Let $\Lambda_0$ be the totality of  increasing sequences $\{\lambda_n\}_{n=1}^{\infty}$
 with $\lim_{n\to\infty} \lambda_n = \infty$ such that the following  $ \lambda $ exists:
$$\lambda :=\lim_{n\to\infty}(\lambda_n-2\pi[\lambda_n/(2\pi)]).$$
For any positive sequence $\{x_n\}_{n=1}^{\infty}$ with $\lim_{n\to\infty} x_n = \infty$, there exists a subsequence $\{\lambda_n\} \in \Lambda_0$ of $\{x_n\}$. We define two positive right-continuous functions $\phi_1(x)$ and  $\phi_2(x)$ on $\R_+$ as
\begin{equation}
\phi_1(x)=e^{-x}(3\pi +1 +\sqrt{2}\sin(x-\frac{\pi}{4}))1_{[0,\infty)}(x) \nonumber
\end{equation}
and
\begin{equation}
\phi_2(x)=\frac{1}{3\pi}1_{[0,2\pi)}(x)+\sum_{n=1}^{\infty}\frac{1}{\pi^3n^2}1_{[2n\pi,2(n+1)\pi)}(x).\nonumber
\end{equation}
Note that the two functions $\phi_1(x)$ and  $\phi_2(x)$ are decreasing on $\R_+$ and $\phi_1(0)\phi_2(0)=1$ and that $\int_0^{\infty}\phi_2(x)dx=1$ and $\int_0^{\infty}\exp(inx)\phi_2(x)dx=0$ for all $n \in \N$. Thus, we can define a distribution $\xi$ on $\R_+$ by using its tail $\bar \xi(x)$  as 
$\bar \xi(x):= \phi_1(x)\phi_2(x)$  on $\R_+$. 
\begin{lem}
We have $\widehat \xi(1)< \infty$ and 
\begin{equation}
\widehat \xi(1 +i)=0.
\end{equation}
\end{lem}

{\it Proof }  Note that $\phi_1(x) \asymp e^{-x}$ and $\phi_2(x) \sim 4\pi^{-1} x^{-2}$, and hence $\bar\xi(x) \asymp e^{-x}x^{-2}$. Thus, $\widehat \xi(1)< \infty$. We have by using integration by parts
\begin{equation}
\widehat \xi(1 +i)= \bar\xi(0) +(1 +i)\int_0^{\infty}e^{(1 +i)x}\bar\xi(x)dx=1-\int_0^{\infty}\phi_2(x)dx=0.\nonumber
\end{equation}  
Thus, the lemma is true. \qed
\begin{lem}
We have $\xi \notin  {\mcal L}(1)$ and hence $\xi \notin  {\mcal S}(1)$. 
\end{lem}
{\it Proof }  For every $\{\lambda_n\}\in \Lambda_0$ and  every $a \in \R$, we have
\begin{equation}
\lim_{n \to \infty}\frac{e^{a}\overline{ \xi}(\lambda_n +a)}{ \overline{ \xi}(\lambda_n )}= \frac{3\pi +1 + \sqrt{2}\sin(\lambda +a -\frac{\pi}{4})}{3\pi +1 + \sqrt{2}\sin(\lambda  -\frac{\pi}{4})},\nonumber
\end{equation}
which is not constant in $a$. Thus, we see that $\xi \notin  {\mcal L}(1)$ and hence $\xi \notin  {\mcal S}(1)$. \qed
\begin{lem}
We have  $\xi^{2*} \in  {\mcal S}(1)$.
\end{lem}
{\it Proof }  Let $g(x):= 1_{[1,\infty)}(x)x^{-2}e^{-x}$ and $A >1.$ Then, we have 
\begin{equation}
\begin{split}
&\lim_{A\to\infty}\limsup_{x \to \infty}\frac{\int_{A}^{x-A}g(x-u)g(u)du}{g(x)}\\
&=\lim_{A\to\infty}\limsup_{x \to \infty}\frac{2\int_{A}^{x/2}g(x-u)g(u)du}{g(x)}\le 8\lim_{A\to\infty}\int_A^{\infty}u^{-2}du=0
\end{split}
\end{equation}
and
\begin{equation}
\lim_{A\to\infty}\lim_{x \to \infty}\frac{g(x-A)g(A)}{g(x)}=\lim_{A\to\infty}A^{-2}=0.
\end{equation}
We see that, for $x > 2A$, 
\begin{equation}
\overline{ \xi^{2*}}(x) =I_1(x) +I_2(x),\nonumber
\end{equation}
where
\begin{equation}
I_1(x) := 2\int_{0-}^{A+}\overline{ \xi}(x-u)\xi(du)\nonumber
\end{equation}
and
\begin{equation}
I_2(x) := \int_{A+}^{(x-A)+}\overline{ \xi}(x-u)\xi(du) +\overline{ \xi}(x-A)\overline{ \xi}(A).\nonumber
\end{equation}
Note that $\overline{ \xi}(x) \asymp g(x)$ and hence $\overline{ \xi}(x) \le c_1 g(x)$ with some $c_1 >0$ for $x >1$ and that $-g'(x) \le 3g(x)$ for $x >1$. By using integration by parts, we obtain that, for $x > 2A$, 
\begin{equation}
\begin{split}
&\int_{A+}^{(x-A)+}\overline{ \xi}(x-u)\xi(du)\\ \nonumber
&\le c_1 \int_{A+}^{(x-A)+}g(x-u)\xi(du)\\ \nonumber
&\le  c_1 g(x-A)\overline{ \xi}(A+)-c_1\int_{A}^{x-A}g'(x-u)\overline{ \xi}(u)du\\ \nonumber
&\le c_1^2g(x-A)g(A)+3c_1^2\int_{A}^{x-A}g(x-u)g(u)du.\nonumber
\end{split}
\end{equation}
Thus, we find that, for $x > 2A$,  
\begin{equation}
I_2(x) \le   2c_1^2g(x-A)g(A)+3c_1^2\int_{A}^{x-A}g(x-u)g(u)du,\nonumber
\end{equation}
and hence by (3.2) and (3.3)
\begin{equation}
\lim_{A\to\infty}\limsup_{x \to \infty}\frac{I_2(x)}{g(x)}=0.
\end{equation}
For every $\{\lambda_n\}\in \Lambda_0$, we have 
\begin{equation}
\lim_{n\to\infty}\frac{I_1(\lambda_n)} {g(\lambda_n)} = 8\pi^{-1}\int_{0-}^{A+}(3\pi +1 + \sqrt{2}\sin(\lambda  -u -\frac{\pi}{4}))e^u\xi(du).\nonumber
\end{equation}
Thus, we see from (3.1) of Lemma 3.1  and (3.4) that, for every $\{\lambda_n\}\in \Lambda_0$,
\begin{equation}
\lim_{n\to\infty}\frac{\overline{ \xi^{2*}}(\lambda_n) }{g(\lambda_n)}=\lim_{A\to\infty}\lim_{n\to\infty}\frac{I_1(\lambda_n)}{g(\lambda_n)}=\frac{8(3\pi+1)}{\pi}\widehat \xi(1),\nonumber
\end{equation}
which is independent of the choice of $\{\lambda_n\}\in \Lambda_0$.  That is, 
\begin{equation}
\overline{ \xi^{2*}}(x) \sim 8(3\pi+1)\pi^{-1}\widehat \xi(1)e^{-x}x^{-2}.\nonumber
\end{equation}
Hence, by (ii) of Lemma 2.1,  we establish that  $\xi^{2*} \in  {\mcal S}(1)$. \qed
 
\medskip

{\it Proof of Theorem 1.1 } The proof is due to Lemmas 3.2 and 3.3. \qed

\medskip

{\it Proof of Corollary  1.1 } The proof of assertion (i)  is due to Theorem 1.1 and (i) of Lemma 2.4. If  the class $ {\mcal S}_{\Delta}$ is closed under convolution roots for some $\Delta$, then so is for every $\Delta$ and thereby the class $ {\mcal S}_{loc}$ is closed under convolution roots. Thus, assertion (ii) is due to assertion (i). The proof of assertion (iii) is due to Lemmas 3.2 and 3.3.  We see from assertion (iii) and (ii) of Lemma 2.4 that  $ {\mcal L}_{loc}\cap {\mcal M}(-\gamma)$ with some $\gamma >0$ is not closed under convolution roots and hence so is $ {\mcal L}_{loc}$. If  the class $ {\mcal L}_{\Delta}$ is closed under convolution roots for some $\Delta$, then so is for every $\Delta$ and thereby the class $ {\mcal L}_{loc}$ is closed under convolution roots. Thus, assertion (v) is due to assertion (iv). \qed

\section{Proofs of Theorem 1.2   and its corollary}
In Sect.\ 4 and 5, let $\mu$ be a distribution on $\R$ satisfying $\bar\mu(x) >0$ for all $x \in \R.$ Let $\{X_j\}_{j=1}^{\infty}$ be IID random variables with distribution $\mu$. We define $J_k(x)$ for $1 \le k \le 3$ and $\ep(A)$ for $A>0$ and $n \ge 2$ as
\begin{equation}
J_1(x):= \int_{-\infty}^{A+}\bar\mu(x-u)\mu^{(n-1)*}(du),\nonumber
\end{equation}
 \begin{equation}
J_2(x):= \int_{A+}^{(x-A)+}\bar\mu(x-u)\mu^{(n-1)*}(du) +\bar\mu(A)\overline{\mu^{(n-1)*}}(x-A),\nonumber
\end{equation}
 \begin{equation}
J_3(x):= \int_{A+}^{(x-A)+}\bar\mu(x-u)\mu(du) +\bar\mu(A)\overline{\mu}(x-A),\nonumber
\end{equation}and
\begin{equation}
\ep(A):= \limsup_{x \to \infty}\frac{J_2(x)+J_3(x)}{\overline{\mu^{n*}}(x)}.\nonumber
\end{equation}
\begin{lem} 
Let $\gamma >0$ and let $n \ge 2$. Then, we have the following.

(i) We have, for $x > nA$,
\begin{equation}
\overline{\mu^{n*}}(x) \le nJ_1(x) + nJ_2(x)
\end{equation}
and
\begin{equation}
nJ_1(x) -2^{-1}n(n-3)J_2(x)-2^{-1}n(n-1)J_3(x)\le \overline{\mu^{n*}}(x). 
\end{equation}

(ii) If $\mu^{n*} \in \mathcal{S}(\gamma)$, then $\lim_{A \to \infty}\ep(A)=0$ and hence we have
\begin{equation}
\lim_{A \to \infty}\liminf_{x \to \infty}\frac{nJ_1(x)}{\overline{\mu^{n*}}(x)}= \lim_{A \to \infty}\limsup_{x \to \infty}\frac{nJ_1(x)}{\overline{\mu^{n*}}(x)}= 1.
\end{equation}
\end{lem}

{\it Proof } \q We have, for $x > nA$, 
\begin{equation}
\begin{split}
\overline{\mu^{n*}}(x)&= P(\sum_{j=1}^n X_j > x)\\ \nonumber
&\le \sum_{k=1}^nP(X_k > A, \sum_{j=1}^n X_j > x)\\ \nonumber
&=  nJ_1(x) + nJ_2(x).
\end{split}
\end{equation}
Thus, (4.1) of assertion (i) is true. On the other hand, we see that, for $x > nA$,
\begin{equation}
\begin{split}
 & P(X_1 > A, X_2 > A,\sum_{j=1}^n X_j > x)\\ \nonumber
& \le P(X_1 > A, X_2 > A,\sum_{j=1}^n X_j > x, \sum_{j=3}^n X_j \ge 0)\\ \nonumber
 & + P(X_1 > A, X_2 > A,\sum_{j=1}^n X_j > x, \sum_{j=3}^n X_j < 0)\\ \nonumber
& \le P(X_1 > A, \sum_{j=2}^n X_j > A, \sum_{j=1}^n X_j > x) + P(X_1 > A, X_2 > A, X_1 +X_2 > x)\\ \nonumber
& = J_2(x) +J_3(x),\\ \nonumber
\end{split}
\end{equation}
with the understanding that $\sum_{j=3}^n X_j=0$ for $n=2$.
Thus,  we have, for $x > nA$,
\begin{equation}
\begin{split}
 & P(\sum_{j=1}^n X_j > x)\\ \nonumber
& \ge \sum_{k=1}^nP(X_k > A, \sum_{j=1}^n X_j > x)-\sum_{1 \le k < l \le n}P(X_k > A, X_l > A,\sum_{j=1}^n X_j > x)\\ \nonumber
& =nJ_1(x) + nJ_2(x)-2^{-1}n(n-1)P(X_1 > A, X_2 > A,\sum_{j=1}^n X_j > x)\\ \nonumber
& \ge nJ_1(x) -2^{-1}n(n-3)J_2(x)-2^{-1}n(n-1)J_3(x).\\ \nonumber
\end{split}
\end{equation}
 Hence, (4.2) of assertion (i) is true. Next, suppose that $\mu^{n*} \in \mathcal{S}(\gamma)$. Let $d:=\mu([0,\infty)).$  We obtain that, for $x > nA$, 
\begin{equation}
\begin{split}
&d^nJ_2(x)+d^{2n-2}J_3(x)\\ \nonumber
= &P(X_1 > A, \sum_{j=2}^n X_j > A,\sum_{j=1}^n X_j > x, X_k \ge0 \mbox{ for } n+1\le k \le 2n )\\ \nonumber
&+P(X_1 > A, X_2 > A,  X_1 +X_2 > x, X_k \ge0 \mbox{ for } 3\le k \le 2n )\\ \nonumber
\le & 2P(X_1 +\sum_{j=n+2}^{2n} X_j> A, \sum_{j=2}^{n+1} X_j > A,\sum_{j=1}^{2n} X_j > x)\\ \nonumber
= & 2\int_{A+}^{(x-A)+}\overline{\mu^{n*}}(x-u)\mu^{n*}(du)+ 2 \overline{\mu^{n*}}(A)\overline{\mu^{n*}}(x-A).\nonumber
\end{split}
\end{equation}
Note from ${\widehat \mu}(\gamma) < \infty$ that 
\begin{equation}
\lim_{A \to \infty}e^{\gamma A }\overline{\mu^{n*}}(A) \le \lim_{A \to \infty}\int_{A+}^{\infty}e^{\gamma x }\mu^{n*}(dx) =0.\nonumber
\end{equation}
Thus, we see from (i) of Lemma 2.2 that
\begin{equation}
\begin{split}
&\lim_{A \to \infty}d^{2n-2}\ep(A)\le \lim_{A \to \infty}\limsup_{x \to \infty}\frac{d^nJ_2(x)+d^{2n-2}J_3(x)}{\overline{\mu^{n*}}(x)}\\ \nonumber
& \le 2\lim_{A \to \infty}\limsup_{x \to \infty}\frac{\int_{A+}^{(x-A)+}{\overline{\mu^{n*}}(x-u)}\mu^{n*}(du)}{\overline{\mu^{n*}}(x)}+2\lim_{A \to \infty}e^{\gamma A }\overline{\mu^{n*}}(A)=0.\nonumber
\end{split}
\end{equation}
Hence, we obtain (4.3) from (4.1) and (4.2). \qed

\medskip

{\it Proof of Theorem 1.2 } \q We continue to use $J_k(x)$ for  $1 \le k \le 3$ and $\ep(A)$ defined above. Let $ n \ge 2$. Define $D^*$ and $D_*$ as
\begin{equation}
D^*:=\limsup_{x \to \infty}\frac{\bar\mu(x)}{\overline{\mu^{n*}}(x)}, \q D_*:=\liminf_{x \to \infty}\frac{\bar\mu(x)}{\overline{\mu^{n*}}(x)}.\nonumber
\end{equation}
Suppose that  $e^{\gamma x}\bar\mu(x) \in {\bf  M}$ and $\mu^{n*} \in \mathcal{S}(\gamma)$. Then, for any $\ep \in (0,1)$ and $A>0$, there is $N> nA$ such that we have, for $x >N$ and $-\infty < u \le A$, 
\begin{equation}
  e^{\gamma (x+A-u)}\bar\mu(x+A-u) \le (1+\ep)e^{\gamma x}\bar\mu(x)
\end{equation}
and, for $x >N$ and $-A \le u \le A$,
\begin{equation}
 (1-\ep) e^{\gamma x}\bar\mu(x) \le e^{\gamma (x-A-u)}\bar\mu(x-A-u). 
\end{equation}
Thus, we see from (4.1) of Lemma 4.1 and (4.4) that
\begin{equation}
\begin{split}
1-n\ep(A)& \le \liminf_{x \to \infty}\frac{nJ_1(x)}{\overline{\mu^{n*}}(x)}\\ \nonumber
&=\liminf_{x \to \infty}\frac{nJ_1(x+A)}{\overline{\mu^{n*}}(x+A)}\\ \nonumber
&=\liminf_{x \to \infty}\frac{n\overline{\mu}(x)}{e^{-\gamma A}\overline{\mu^{n*}}(x)}\int_{-\infty}^{A+}\frac{\overline{\mu}(x+A-u)}{\overline{\mu}(x)}\mu^{(n-1)*}(du)\\ \nonumber
&\le nD_*\int_{-\infty}^{A+}e^{\gamma u}\mu^{(n-1)*}(du). \nonumber
\end{split}
\end{equation}
On the other hand, we obtain from  (4.2) of Lemma 4.1 and (4.5) that 
\begin{equation}
1 +2^{-1}n(n-1) \ep(A) \ge \limsup_{x \to \infty}\frac{nJ_1(x-A)}{\overline{\mu^{n*}}(x-A)}\ge nD^*\int_{(-A)+}^{A+}e^{\gamma u}\mu^{(n-1)*}(du).\nonumber
\end{equation}
As $A \to \infty$, we have, by (ii) of Lemma 4.1,
\begin{equation}
D^*= D_*= n^{-1}{\widehat \mu}(\gamma)^{1-n}.\nonumber
\end{equation}
Hence, we establish that 
\begin{equation}
\bar\mu(x)\sim n^{-1}{\widehat \mu}(\gamma)^{1-n}\overline{\mu^{n*}}(x).\nonumber
\end{equation}
Thus, we conclude from (ii) of  Lemma 2.2 that $\mu \in \mathcal{S}(\gamma)$. \qed

\medskip

{\it Proof of Corollary 1.2 } The proof is due to Theorem 1.2 and Lemma 2.3. \qed

\medskip

\section{Proofs of Theorem 1.3  and its corollary}
Let $\Lambda_1$ be the totality of  increasing sequences $\{\lambda_k\}_{k=1}^{\infty}$
 with $\lim_{k\to\infty} \lambda_k = \infty$ such that, 
 for every  $x\in \R$,  the following  $  m(x;\{\lambda_k\} )$ exists and is finite:
\begin{equation}
m(x;\{\lambda_k\}):=\lim_{k\to\infty}\frac{\overline{\mu}(\lambda_k+x)}{\overline{\mu^{n*}}(\lambda_k)}.
\end{equation}
The idea of the use of the function $ m(x;\{\lambda_k\})$ goes back to Teugels \cite{te75} and is extensively employed in Watanabe and Yamamuro \cite{wy10a}.
\begin{lem} Assume that $\mu^{n*} \in \mathcal{S}(\gamma)$ for $n \ge 2$. Define $d:=\mu([0,\infty))$. For any sequence $\{x_k\}_{k=1}^{\infty}$ with $\lim_{k\to\infty} x_k = \infty$, there exists a subsequence  $\{\lambda_k\} \in \Lambda_1$ of $\{x_k\}$. Moreover,  $m(x;\{\lambda_k\} )$ is decreasing and finite, and  we have
\begin{equation}
M(x;\{\lambda_k\}):= e^{\gamma x} m(x;\{\lambda_k\} ) \le d^{1-n}. 
\end{equation}

\end{lem}

{\it Proof } Let 
$$T_k(y):=\frac{\overline{\mu}(x_k+y)}{\overline{\mu^{n*}}(x_k)}.$$
Then,  we see that $T_k(y)$ is decreasing and
\begin{equation}
\sup_{y \in [-A,A]}T_k(y)\le \frac{\overline{\mu}(x_k-A)}{\overline{\mu^{n*}}(x_k)} \le d^{1-n}\frac{\overline{\mu^{n*}}(x_k-A)}{\overline{\mu^{n*}}(x_k)} \nonumber
\end{equation}
and
\begin{equation}
\lim_{k\to\infty}\frac{\overline{\mu^{n*}}(x_k-A)}{\overline{\mu^{n*}}(x_k)} = e^{\gamma A}.\nonumber
\end{equation}
Thus, $T_k(y)$ is uniformly bounded on all finite intervals. By the selection principle (See Chap. VIII of Feller \cite{feller1}),  there exists an increasing subsequence $\{\lambda_k\} \in \Lambda_1$ of $\{x_k\}$. 
Since $T_k(x)$ is decreasing, so is  $m(x;\{\lambda_k\} )$. Moreover, 
\begin{equation}
M(x;\{\lambda_k\})\le d^{1-n} \lim_{k\to\infty}\frac{ e^{\gamma x} \overline{\mu^{n*}}(\lambda_k+x)}{\overline{\mu^{n*}}(\lambda_k)}= d^{1-n}. \nonumber
\end{equation}
Thus, the lemma is true.  \qed

\medskip

\begin{lem} Let $\gamma >0.$ Define  $f_c(x): =c_1^{-1}e^{-\gamma x}(1-c^{-1}x)1_{[0,c)}(x)$ for $ c >0$  with $c_1:= \int_o^ce^{-\gamma x}(1-c^{-1}x)dx$ and define an absolutely continuous distribution $\mu_c$ as $ \mu_c:= (f_c(x)dx)*\mu $ for $ c >0$. Then, we have the following.

(i) Fix $c >0$.   We have $\widehat\mu(\gamma ) < \infty$ and $\widehat\mu(\gamma + iz) \neq 0$ for every $z \in \R$ if and only if  $\widehat\mu_c(\gamma ) < \infty$ and $\widehat\mu_c(\gamma + iz) \neq 0$ for every $z \in \R$.

(ii) Let $n \in \N$. We have  $\mu^{n*} \in \mathcal{S}(\gamma)$ if and only if $\mu_c^{n*} \in \mathcal{S}(\gamma)$  for every $c >0$.
 \end{lem}

{\it Proof } Note that $\int_o^ce^{\gamma x}f_c(x)dx < \infty$ and $\int_o^ce^{\gamma x + izx}f_c(x)dx \neq 0$ for every $z \in \R$. Thus, the proof of assertion (i) is obvious. Next, we prove assertion (ii) for $n=1$. 
Suppose that $\mu \in \mathcal{L}(\gamma)$. Then, we  have, for every $c >0$,
\begin{equation}
\lim_{x\to\infty}\frac{\overline{\mu_c}(x)}{\overline{\mu}(x)} = \lim_{x\to\infty}\int_0^c\frac{\overline{\mu}(x-u)}{\overline{\mu}(x)} f_c(u)du=\int_0^ce^{\gamma u}f_c(u)du. 
\end{equation}
Thus, by (ii) of Lemma 2.2, if $\mu \in \mathcal{S}(\gamma)$, then $\mu_c \in \mathcal{S}(\gamma)$  for every $c >0$. Conversely, suppose that $\mu_c \in \mathcal{S}(\gamma)$ for every $c >0$.
 Let $X$ be a random variable with distribution $\mu$ and let $Y_c$ be a random variable with distribution $f_c(x)dx$ independent of $X$. Then, we have 
\begin{equation}
P(x +c< X +Y_c) \le P( x < X) \le P(x < X +Y_c). \nonumber
\end{equation}
Thus, we see that, for every $a \in \R$,
\begin{equation}
\begin{split}
e^{-\gamma (a +c)} &= \lim_{x\to\infty} \frac{P(x +a+c< X +Y_c)}{P( x < X +Y_c)}\\ \nonumber
&\le \liminf_{x\to\infty}\frac{P(x+a < X )}{P( x < X)} \\ \nonumber
&\le\limsup_{x\to\infty}\frac{P(x+a < X )}{P( x < X)} \\ \nonumber
&\le\lim_{x\to\infty}\frac{P(x+a < X +Y_c)}{P( x+c < X+Y_c)} =e^{-\gamma (a -c)}.\\ \nonumber
\end{split}
\end{equation}
As $c \downarrow 0$, we find that, for every $a \in \R$,
\begin{equation}
\lim_{x\to\infty}\frac{P(x+a < X )}{P( x < X)}= e^{-\gamma a}.\nonumber
\end{equation}
Hence,  $\mu \in \mathcal{L}(\gamma)$ and thereby we establish from (5.3) and (ii) of Lemma 2.2 that $\mu \in \mathcal{S}(\gamma)$.  The proof of assertion (ii) for $ n \ge 2$ is similar and omitted. \qed

\medskip

{\it Proof of Theorem 1.3 } Suppose that $\mu^{n*} \in \mathcal{S}(\gamma)$ for $n \ge 2$ and that  $\widehat\mu(\gamma ) < \infty$ and $\widehat\mu(\gamma + iz) \neq 0$ for every $z \in \R$. Let $\mu_c$ be the absolutely continuous distribution defined in Lemma 5.2. Then, by  Lemma 5.2, we see that  $\mu_c^{n*} \in \mathcal{S}(\gamma)$ for every $c >0$ and that $\widehat\mu_c(\gamma ) < \infty$ and $\widehat\mu_c(\gamma + iz) \neq 0$ for every $z \in \R$ and  every $c >0$. By using $\mu_c$  instead of $\mu$ in the definitions of (5.1) and (5.2), we replace the class $\Lambda_1$ and the functions $  m(x;\{\lambda_k\} )$ and  $  M(x;\{\lambda_k\} )$ by $\Lambda_{1,c}$,  $  m_c(x;\{\lambda_k\} )$, and  $  M_c(x;\{\lambda_k\} )$, respectively.  Let $A >0$ and $a \in \R$.  Define, for  $\{\lambda_k\} \in \Lambda_{1,c}$, 
$$U_k(y):=\frac{\overline{\mu_c}(\lambda_k+y)}{\overline{\mu_c^{n*}}(\lambda_k)}.$$
  Since $U_k(a-u)$ is bounded on $(-\infty, A]$, we obtain from the dominated convergence theorem and (4.3) of Lemma 4.1 that
\begin{equation}
\begin{split}
e^{-\gamma a}=\lim_{k\to\infty}\frac{\overline{\mu_c^{n*}}(\lambda_k+a)}{\overline{\mu_c^{n*}}(\lambda_k)} &= \lim_{A \to \infty}n \int_{-\infty}^{A+}\lim_{k\to\infty}U_k(a-u)\mu_c^{(n-1)*}(du)\\ \nonumber
&=n \int_{-\infty}^{\infty}m_c(a-u;\{\lambda_k\} )\mu_c^{(n-1)*}(du).
\end{split}
\end{equation}
Thus, we find that, for every $a \in \R$,
\begin{equation}
1= n \int_{-\infty}^{\infty} M_c(a-u;\{\lambda_k\} )e^{\gamma u}\mu_c^{(n-1)*}(du).
\end{equation}
Hence, we have, for every $a, b \in \R$,
\begin{equation}
  \int_{-\infty}^{\infty} (M_c(a+b-u;\{\lambda_k\} )-M_c(b-u;\{\lambda_k\} ))e^{\gamma u}\mu_c^{(n-1)*}(du)=0.\nonumber
\end{equation}
Note  that, by (5.2),  $M_c(x;\{\lambda_k\} )$ is bounded and that, for every $z \in \R$,
$$\int_{-\infty}^{\infty}e^{izx}e^{\gamma x}\mu_c^{(n-1)*}(dx)=\widehat\mu_c(\gamma +iz)^{n-1}\ne 0.$$ 
It follows from Lemma 2.5 that, for every $a \in \R$,
\begin{equation}
   M_c(a+b;\{\lambda_k\} )-M_c(b;\{\lambda_k\} )=0 \mbox{ for a.e. } b \in \R.\nonumber
\end{equation}
Since the function $  m_c(x;\{\lambda_k\} )$ is decreasing, the functions $  M_c(x-;\{\lambda_k\} )$ and $  M_c(x+;\{\lambda_k\} )$ exist for all $x \in \R$.
Taking $b_n = b_n(a) \uparrow  0$ and  $b_n = b_n(a)\downarrow 0$, we have, for every $a \in \R$, 
\begin{equation}
   M_c(a-;\{\lambda_k\} )=M_c(0-;\{\lambda_k\} )\mbox{ and } M_c(a+;\{\lambda_k\} )=M_c(0+;\{\lambda_k\} ).
\end{equation}
Then, taking $a \downarrow 0$ in the first equality of (5.5), we have 
$$C:=M_c(0+;\{\lambda_k\} )=M_c(0-;\{\lambda_k\} ).$$
Thus, we obtain from (5.5) that,  for every $a \in \R$,
$$M_c(a;\{\lambda_k\} )=M_c(a-;\{\lambda_k\} )=M_c(a+;\{\lambda_k\} )=C.$$
Therefore, we see from (5.4) that 
\begin{equation}
C=\lim_{k\to\infty}\frac{\overline{\mu_c}(\lambda_k)}{\overline{\mu_c^{n*}}(\lambda_k)}= n^{-1} \widehat\mu_c(\gamma)^{1-n},\nonumber
\end{equation}
which is independent of the choice of $\{\lambda_k\} \in \Lambda_{1,c}$. Thus, we find from Lemma 5.1 that 
$$\overline{\mu_c}(x)\sim n^{-1}\widehat\mu_c(\gamma)^{1-n}\overline{\mu_c^{n*}}(x)$$
and, by (ii) of Lemmas 2.2,  $\mu_c \in \mathcal{S}(\gamma)$ for every $c >0$. Thus, we conclude from (ii) of Lemma 5.2 for $n=1$ that  $\mu \in \mathcal{S}(\gamma)$. \qed
\medskip

{\it Proof of Corollary 1.3 } The proof is due to Theorem 1.3 and Lemma 2.3. \qed


\medskip


\begin{thebibliography}{22}

 

\bibitem{afk03} Asmussen, S.,  Foss, S.,  Korshunov, D. : Asymptotics for sums of random variables with local subexponential behaviour.  J. Theoret. Probab.  {\bf 16},     489-518 (2003) 

\bibitem{bi87}   Bingham,  N. H.,  Goldie,  C. M.,  Teugels,  J. L. :  Regular Variation.  Cambridge University Press, Cambridge. (1987)


\bibitem{cnw} Chover, J.; Ney, P., Wainger, S. : Functions of probability measures. J. Analyse Math. {\bf 26},  255-302 (1973) 

\bibitem{cnw1} Chover, J.; Ney, P.,  Wainger, S. : Degeneracy properties of subcritical branching processes. Ann. Probab. {\bf 1},  663-673 (1973)

\bibitem{eg80} Embrechts, P., Goldie, C. M. :  On closure and factorization properties of subexponential and related distributions.  J. Austral. Math. Soc. Ser. A  {\bf 29},  243-256 (1980)

\bibitem{eg82}  Embrechts, P.,  Goldie, C. M. :  On convolution tails.  Stochastic Process. Appl.  {\bf 13},  263-278 (1982)

\bibitem{egv79} Embrechts, P.,  Goldie,  C. M.,  Veraverbeke, N. :  Subexponentiality and infinite divisibility.  Z. Wahrsch. Verw. Gebiete  {\bf 49}, 335-347 (1979)

\bibitem{feller1} Feller, W. :  An Introduction to Probability Theory and its Applications. Vol. II, 2nd ed. Wiley, New York. (1971) 
 
\bibitem{fkz13} Foss, S.,  Korshunov, D., Zachary, S. : An introduction to heavy-tailed and subexponential distributions. Second edition. Springer Series in Operations Research and Financial Engineering. Springer, New York. (2013)

\bibitem{kl89}  Kl\"uppelberg, C. :   Subexponential distributions and characterizations of related classes.  Probab. Theory Related Fields  {\bf 82},  259-269 (1989)

\bibitem{ko04} Korevaar, J. : Tauberian theory. A century of developments. Grundlehren der Mathematischen Wissenschaften  329. Springer-Verlag, Berlin. (2004)

\bibitem{pa04} Pakes, A. G. :  Convolution equivalence and infinite divisibility.  J. Appl. Probab.  {\bf 41},  407-424 (2004)

\bibitem{pa07} Pakes, A. G.:  Convolution equivalence and infinite divisibility: Corrections and corollaries.  J. Appl. Probab. {\bf 44},  295-305  (2007)

\bibitem{sato} Sato, K. :  L\'{e}vy Processes and Infinitely Divisible Distributions. Cambridge University Press, Cambridge. (1999)

\bibitem{shw05} Shimura, T.,   Watanabe,  T. :  Infinite divisibility and generalized subexponentiality.  Bernoulli  {\bf 11},  445-469 (2005)

\bibitem{shw} Shimura, T.,  Watanabe, T. :   On the convolution roots in the convolution-equivalent class. The Institute of Statistical Mathematics Cooperative Research Report 175 pp1-15 (2005)

\bibitem{te75} Teugels, J. L. :  The class of subexponential distributions.  Ann. Probab.  {\bf 3},  1000-1011 (1975)

\bibitem{wa08} Watanabe, T. : Convolution equivalence and distributions of random sums. Probab. Theory Related Fields {\bf 142},  367-397 (2008)   

\bibitem{wy09} Watanabe, T., Yamamuro, K. : Local subexponentiality of infinitely divisible distributions. J. Math-for-Ind. {\bf 1},  81-90 (2009)

\bibitem{wy10} Watanabe, T., Yamamuro, K. :  Local subexponentiality and self-decomposability. J. Theoret. Probab. {\bf23},  1039-1067 (2010) 


\bibitem{wy10a} Watanabe, T., Yamamuro, K. : Ratio of the tail of an infinitely divisible distribution on the line to that of its L\'evy measure. Electron. J. Probab. {\bf 15},  44-74 (2010) 


\bibitem{w32} Wiener, N. : Tauberian theorems. Ann. of Math. (2) {\bf 33},   1-100 (1932)

\end{thebibliography}
\end{document}